\documentclass[12pt]{amsart}
\newcommand{\bC}{{\mathbb C}}

\newcommand{\bQ}{{\mathbb Q}}

\newcommand{\cI}{{\mathcal I}}

\newcommand{\cO}{{\mathcal O}}

\newcommand{\cZ}{{\mathcal Z}}

\newcommand{\tA}{\widetilde{A}}

\newcommand{\tD}{\widetilde{D}}

\newcommand{\tZ}{\widetilde{Z}}
\newcommand{\ra}{\rightarrow}
\newcommand{\lra}{\longrightarrow}

\newcommand{\inj}{\hookrightarrow}

\newtheorem{proposition}{Proposition}[section]
\newtheorem{lemma}[proposition]{Lemma}

\newtheorem{theorem}[proposition]{Theorem}
\newtheorem{theoremi}{Theorem}
\newtheorem{conjecture}[theoremi]{Conjecture}
\newtheorem{definition}[proposition]{Definition}

\newtheorem{remark}[proposition]{Remark}

\numberwithin{equation}{section}

\begin{document}

\title[]{An inductive approach to the Hodge conjecture for abelian varieties}

\author{E. Izadi}

\address{Department of Mathematics, Boyd
Graduate Studies Research Center, University of Georgia, Athens, GA
30602-7403, USA}

\email{izadi@math.uga.edu}

\thanks{This material is based upon work partially supported by the
National Security Agency (NSA) and the
National Science Foundation (NSF). Any opinions,
findings and conclusions or recomendations expressed in this material
are those of the author and do not necessarily reflect the views of
NSF or NSA}


\maketitle

\section*{Introduction}

Let $X$ be a smooth complex projective variety of dimension $g$. A
Hodge class of degree $2d$ on $X$ is, by definition, an element of
$H^{2d} (X, \bQ )\cap H^{ d,d} (X)$. The cohomology class of an
algebraic subvariety of codimension $d$ of $X$ is a Hodge class of
degree $2d$. The original Hodge conjecture states that any Hodge class
on $X$ is algebraic, i.e., a $\bQ$-linear combination of classes of
algebraic subvarieties of $X$. Lefschetz' Theorem says that Hodge
classes of degree $2$ are always algebraic.

The classical Hodge conjecture has been generalized by Grothendieck as
follows (see Steenbrink \cite{steenbrink87} page 166). To fix some
notation, we will always designate a Hodge structure by its rational
vector space $V$, the splitting $V\otimes\bC =\oplus_{ p+q =m} V^{p,q}$
being implicit. We say that $V$ is effective if $V^{ p,q} =0$ when
either $p$ or $q$ is negative. Recall that the level of a Hodge
structure is the integer
\[
Max\{ |p-q| : V^{ p,q}\neq 0\}
\]
so that to say that $V$ has level $l$ means that after a Tate twist,
$V$ will be an effective Hodge structure of weight $l$ with
non-vanishing $(l,0)$ component.

Below $m$ and $p$ are two positive integers.

\begin{conjecture}
{\bf GHC(X,m,p):} For every $\bQ$-Hodge substructure $V$ of $H^m
(X,\bQ )$ with level $\leq m- 2p$, there exists a subvariety $Z$ of
$X$ of pure codimension $p$ such that $V\subset Ker\{ H^m (X,\bQ )\ra H^m (
X\setminus Z,\bQ )\}$.
\end{conjecture}
Given $V$ and $Z$ as in the conjecture, we say that $V$ is supported
on $Z$ or that $Z$ supports $V$. Letting $\tZ\ra Z$ be a generically
finite morphism from a nonsingular variety $\tZ$, we have the Gysin
map $H^{ m-2p} (\tZ ,\bQ
)\ra H^m (X ,\bQ)$. The above conjecture is equivalent to
\begin{conjecture}
{\bf GHC(X,m,p):} For every $\bQ$-Hodge substructure $V$ of $H^m (X,\bQ )$
with level $\leq m- 2p$, there exists a subvariety $Z$ of $X$ of pure
codimension $p$ such that $V$ is contained in the image of the Gysin map
$H^{ m-2p} (\tZ ,\bQ )\ra H^m (X ,\bQ)$.
\end{conjecture}

For a smooth ample divisor $Y$ in $X$ we prove below that certain
cases of the Hodge conjectures are equivalent to or a consequence of
those for $Y$. Under some additional assumptions made precise below we
give a simple inductive procedure which relates the Hodge conjectures
for $X$ to those for $Y$. We then apply this to abelian varieties.

\section{The induction}

The following results were inspired by a theorem of Grothendieck and
Steenbrink (Theorem 2.2 in \cite{steenbrink87} page 167). Assume that
$Y$ is a smooth ample divisor in $X$. Suppose $X$ has dimension $g$,
$n$ is a positive integer and $p$ is a positive integer less than or
equal to $\frac{n}{2}$.
\begin{lemma}\label{lemfund}
If $n\leq g-2$, then
\[
GHC(Y, n, p)\Longleftrightarrow GHC (X, n, p)\; ,
\]
and,
\[
GHC (Y,g-1, p)\Longrightarrow GHC (X , g-1 ,p)\; .
\]
\end{lemma}
\begin{proof} Suppose $n\leq g-1$ and let us prove the implications
$\Rightarrow$. Let $V\subset H^n (X ,\bQ)$ be a Hodge substructure of
level $\leq n-2p$. Then, since by the weak Lefschetz theorem, $H^n (X
,\bQ)\inj H^n (Y ,\bQ )$, we can consider $V$ to be a Hodge
substructure of $H^n (Y ,\bQ)$. So by assumption there is a subvariety
$Z$ of $Y$, of pure codimension $p$, such that $V\subset Ker\{ H^n (Y
,\bQ)\ra H^n (Y\setminus Z ,\bQ )\}$. Let $W$ be a pure codimension $p$
subvariety of $X$ containing $Z$ such that $Y\cap W$
has pure codimension $p+1$ in $X$. Then $H^n (Y ,\bQ)\ra H^n (Y\setminus
W ,\bQ)$ factors through $H^n (Y ,\bQ)\ra H^n (Y\setminus Z ,\bQ )$
and hence its kernel contains $V$. Consider the following commutative
diagram of pull-back maps on cohomology
\[\begin{array}{ccc}
H^n (X ,\bQ ) &\inj & H^n (Y ,\bQ)\\
\downarrow & &\downarrow\\
H^n (X\setminus W ,\bQ ) &\lra & H^n (Y\setminus W ,\bQ ).
\end{array}
\]

\begin{lemma}\label{leminj}
The map $H^n (X\setminus W ,\bQ ) \ra  H^n (Y\setminus W ,\bQ )$ in
the above diagram is injective.
\end{lemma}
\begin{proof}
This follows from \cite{goreskymacpherson88}, Theorem
pages 150-151: In a projective embedding of $X$ given by some
multiple of $Y$, let $H$ be a hyperplane cutting a multiple of $Y$ on
$X$. Then, since $Y$ is smooth and $Y\cap W$ has codimension at least
$2$ in $X$, $Y\setminus W$ is a deformation retract of a
$\delta$-neighborhood of $H\cap (X\setminus W)$ in $X\setminus W$.
\end{proof}

Therefore
$V\subset Ker\{ H^n (X ,\bQ)\ra H^n (X\setminus W ,\bQ )\}$.

Now suppose $n\leq g-2$ and let us prove the implication
$\Leftarrow$. In this case, by the weak Lefschetz theorem, the
pull-back map $H^n (X ,\bQ )\ra H^n (Y ,\bQ )$ is an isomorphism. Let
$V\subset H^n (Y ,\bQ)$ be a Hodge substructure of level $\leq
n-2p$. Let $W\subset X$ be a pure codimension $p$ subvariety such that
$V\subset Ker\{ H^n (X ,\bQ)\ra H^n (X\setminus W ,\bQ )\}$. Then a
commutative diagram such as the above shows that $Z :=Y\cap W$ has the
property $V\subset Ker\{ H^n (Y ,\bQ)\ra H^n (Y\setminus Z ,\bQ
)\}$. If $Z$ has pure codimension $p$ in $Y$ then we are done. If not, let
$Y'$ be a smooth ample divisor in $X$ such that $Y\cap Y'$ is smooth
and $Y'$ intersects $W$ transversely. Then every component of $Z' :=
Y\cap Y'\cap W$ has codimension at least $p$ in $Y$. Replace the
components of codimension $>p$ of $Z'$ by irreducible subvarieties of
codimension $p$ containing them. Then, since $n\leq g-2$, as above we
have the commutative diagram with top horizontal injective map
\[
\begin{array}{ccc}
H^n (X ,\bQ) &\inj & H^n (Y\cap Y' ,\bQ ) \\
\downarrow & & \downarrow \\
H^n (X\setminus W ,\bQ) &\lra & H^n (Y\cap Y'\setminus Z' ,\bQ )
\end{array}
\]
which shows
\[
V\subset Ker\{ H^n (Y\cap Y' ,\bQ)\lra H^n (Y\cap Y'\setminus Z',\bQ ).
\]
Now using the diagram
\[
\begin{array}{ccc}
H^n (Y ,\bQ) &\inj & H^n (Y\cap Y' ,\bQ ) \\
\downarrow & & \downarrow \\
H^n (Y\setminus Z' ,\bQ) &\inj & H^n (Y\cap Y'\setminus Z' ,\bQ )
\end{array}
\]
where the injectivity of the bottom horizontal map is proved as in
Lemma \ref{leminj}, we obtain
\[
V\subset Ker\{ H^n (Y ,\bQ)\lra H^n (Y\setminus Z',\bQ )\}
\]
and we are done.
\end{proof}

The primitive part $K(Y, \bQ)$ of the cohomology of $Y$ can be
defined as the kernel of the Gysin map $H^{g-1}(Y, \bQ) \lra H^{g+1}(X,
\bQ)$ obtained by Poincar\'e Duality from push-forward on homology.

\begin{definition}
We shall say that $GHC' (Y, p)$ holds if every Hodge substructure of level
$\leq g-1-2p$ of $K (Y ,\bQ)$ is supported on a pure codimension $p$
subvariety of $Y$.
\end{definition}

Complementary to the lemma is the following result.

\begin{proposition}\label{propGHT}
We have
\[
GHC (Y ,g-1 ,p)\Longleftrightarrow GHC (X, g-1 ,p)\hbox{ and }GHC' (Y ,
p)\; .
\]
\end{proposition}

\begin{proof} Similar to the proof of Lemma \ref{lemfund} above.\end{proof}

Now assume there is a smooth and complete curve $C$, a smooth variety $Y'$ of
dimension $g-1$ and a surjective morphism $C\times Y'\ra X$ such that
there exists $o\in C$ such that the image of $\{ o\}\times Y'$ in $X$
is $Y$. Our main examples for varieties with such a property are symmetric
powers of curves and abelian varieties:

For a curve $T$ of any genus, the symmetric power $T^{ (g-1) }$ embeds
as a smooth ample divisor in $X = T^{ (g) }$ via addition of a point of
$T$. We can then take $Y= Y' = T^{ (g-1) }$ and $C=T$.

For an abelian variety $A=X$, we can take $Y=Y'$ to be any smooth
ample divisor in $A$ and $C$ to be a smooth curve with a generically
injective map $C\ra A$ whose image generates $A$ as a group.

With the above assumptions on $X, Y, Y'$ and $C$ we have

\begin{theorem}\label{thmGHT}
Suppose $GHC(Y' ,g-1, p-1)$ and $GHC (Y' ,g-2, p-1)$ hold and, for any
subvariety $Z$ of $Y'$ of pure codimension $p-1$, $GHC (C\times\tZ, g-2p+2
,1)$ holds, then $GHC( X, g, p)$ holds. Moreover, if $g = 2p+2$, then
$GHC (X ,g, p+1)$ also holds.
\end{theorem}

\begin{remark} Note that the theorem is only interesting if $p\geq
2$. For $g=4$, the theorem is especially interesting since the
hypothesis $GHC (C\times\tZ, g-2p+2 ,1)$ for all $Z\subset Y'$ of pure
codimension $p-1$ is simply Lefschetz' Theorem for $(1,1)$-classes and
is automatically true.
\end{remark}

\begin{remark}

For $g$ even and $p=\frac{g}{2}$ this is reminiscent of the celebrated
approach of Poincar\'e-Lefschetz-Griffiths to the classical Hodge
conjecture via the theory of normal functions (see for instance
\cite{zucker76}, and \cite{lewis99} Chapters $12$ and $14$). Note,
however, that the theorem on normal functions requires very ample
divisors, whereas this approach can be applied to small ample
divisors. The geometry and cohomology of smaller divisors is simpler
and more manageable.

\end{remark}

\begin{proof} The pull-back
\[\begin{split}
H^g ( X,\bQ )\lra H^g (C\times Y' ,\bQ )\cong \\
H^0 (C ,\bQ )\otimes_{\bQ}
H^g (Y' ,\bQ )\oplus H^1 (C ,\bQ )\otimes_{\bQ} H^{g-1} (Y' ,\bQ )\oplus
H^2 (C ,\bQ )\otimes H^{ g-2 } (Y' ,\bQ)
\end{split}
\]
is injective. Using cup product to identify $H^g (X ,\bQ)$ and $H^g
(C\times Y' ,\bQ)$ with their
duals, the transpose of pull-back is Gysin push-forward and is
surjective. Let $V$ be a Hodge substructure of level $\leq g-2p$ of
$H^g ( X,\bQ )$. Let $V_0, V_1$ and $V_2$ be the images of $V$ by the
compositions
\[
V\lra H^g ( X,\bQ ) \lra H^g ( C\times Y',\bQ )\lra H^0 (C ,\bQ
)\otimes_{\bQ} H^g ( Y' ,\bQ ),
\]
\[
V\lra H^g ( X,\bQ ) \lra H^g ( C\times Y',\bQ )\lra H^1 (C ,\bQ )\otimes_{\bQ}
H^{g -1} ( Y' ,\bQ )
\]
and
\[
V\lra H^g ( X,\bQ ) \lra H^g ( C\times Y',\bQ )\lra H^2 (C ,\bQ )\otimes
H^{ g-2 } ( Y' ,\bQ)
\]
respectively. If each of $V_0$, $V_1$ and $V_2$ is supported on a
subvariety of codimension $p$ of $C\times Y'$, then $V$ is supported on
the union of the images of these subvarieties in $X$ and the Hodge
conjecture will follow for $V$. First consider
\[
V_2\subset H^2 (C ,\bQ )\otimes H^{ g-2 } ( Y' ,\bQ)\cong H^{ g-2 } ( Y'
,\bQ).
\]
Since $g-2p = g-2 -2(p-1)$, by assumption there is a subvariety $Z_2$
of pure codimension $p-1$ of $Y'$ which supports $V_2$. Let $t$ be a
general point of $C$. Then $ \{ t\}\times Z_2\subset C\times Y'$ has
pure codimension $p$ and supports $V_2$.

Next consider
\[
V_0\subset H^0 (C ,\bQ )\otimes H^g ( Y' ,\bQ)\cong H^g
( Y' ,\bQ).
\]
Choose an ample divisor $D$ on $Y'$. Cup product with the class of $D$
induces the isomorphism $\cup [D] : H^{ g-2 } (Y'
,\bQ)\stackrel{\sim}{\ra} H^g (Y' ,\bQ)$. Let $Z_0'$ be a subvariety
of pure codimension $p-1$ of $Y'$ supporting the inverse image of
$V_0$ under the isomorphism $\cup [D]$. For $m$ large enough, a
sufficiently general divisor $E\in |\cO_{Y' } (mD) |$ is transverse to
every element of the image of $H^{ g- 2p }(\tZ_0' ,\bQ)$. Hence
$V_0\subset H^g (Y' ,\bQ)$ is supported on $Z_0 := E\cap Z_0'$ which
has pure codimension $p$ in $Y'$. The map $H^{g- 2p } (C\times\tZ_0 )\ra
H^g (C\times Y')$ preserves K\"unneth components hence $V_0\subset H^g
(C\times Y')$ is supported on $C\times Z_0$.

Finally we consider
\[
V_1\subset H^1 (C ,\bQ )\otimes_{\bQ} H^{ g-1 } (Y',\bQ).
\]
Using the intersection pairing to identify $H^1 (C ,\bQ )$ with its
dual, we obtain the map
\[
V_1\otimes H^1 (C ,\bQ )\lra H^{g -1}(Y' ,\bQ )
\]
whose image $W$ is a Hodge substructure of level $\leq g-2p+1 = g-1-
2(p-1)$. So $W$ is supported on a subvariety $Z_1$ of pure codimension
$p-1$ of $Y'$. Then $V_1\subset H^1 (C ,\bQ )\otimes W$ is supported
on $C\times Z_1$. If $GHC(C\times \tZ_1, g-2p+1, 1)$ holds, then $V_1$
is supported on a subvariety of pure codimension $1$ of $C\times\tZ_1$
whose image in $C\times Y'$ has pure codimension $p$. The image of this
subvariety in $X$ is the solution to the Hodge conjecture in this
case.

The last assertion now follows from Lemma 2.1 on page 167 of
\cite{steenbrink87}.
\end{proof}

\section{Abelian varieties with an action of an imaginary quadratic field}

Since the Hodge conjecture is known for degree two classes, it is
also known for Hodge classes that are polynomial combinations of
degree $2$ Hodge classes. Hodge classes which are {\em not}
combinations of degree $2$ classes are usually called {\em
exceptional}.

Perhaps the first example of an abelian variety with an exceptional
Hodge class was that of an abelian fourfold with complex
multiplication given by Mumford \cite{mumford66}. Weil
\cite{weil79III1} observed that the presence of the exceptional Hodge
classes was not due to $A$ being of CM-type but to the fact that the
$\bQ$-endomorphism ring of $A$ contained an imaginary quadratic field
$K$, stable under all Rosati involutions, whose action on the tangent
space of $A$ had eigenspaces of equal dimension, say $n$. Note that
the dimension of the abelian variety is then $g=2n$. These are now
called abelian varieties of Weil type. Weil proved that $\wedge^{2n
}_K H^1 (A ,\bQ) \subset H^{2n} (A, \bQ)\cap H^{ n,n }(A)$, i.e.,
consists of Hodge classes (see \cite{weil79III1} and
\cite{geemen94} Lemma 5.2 page 238). The elements of $W^0$ are called
Weil classes. We observe:

\begin{lemma}
For each $m\leq g$, the
exterior power $W^m :=\wedge^{2n -m}_K H^1 (A ,\bQ)$ over $K$ is a Hodge
substructure of $H^{2n-m} (A, \bQ)$ of level $m$.
\end{lemma}

\begin{proof}
The proof of this is entirely similar to the proof of Lemma 5.2 on
page 238 of \cite{geemen94}.
\end{proof}

So the cohomology of an abelian variety of Weil type contains many
extra Hodge substructures (as compared to a generic abelian variety).

More generally, we make the following

\begin{definition}
An abelian variety $A$ is of Weil type $k\geq 0$ if the action of $K\inj
End_0 (A)$ is stable under all Rosati involutions and the eigenspaces
of the action of $K$ on $T_0A$ have dimension $n$ and $n+k$.
\end{definition}

Then a similar proof gives

\begin{lemma}
Suppose $A$ is of Weil type $k$. For each $m\leq g$, the
exterior power $W^m :=\wedge^{g -m}_K H^1 (A ,\bQ)$ over $K$ is a Hodge
substructure of $H^{g-m} (A, \bQ)$ of level $m+k$.
\end{lemma}

Weil also showed that any abelian variety of Weil type $k=0$ is a member of
a family of dimension $n^2$ of such abelian varieties. By a general
abelian variety of Weil type we mean a general member of such a
family. Weil proved that Weil classes on general abelian varieties of
Weil type are exceptional
\cite{weil79III1}.

In the case where $k\neq 0$, one can similarly construct families of
dimension $n(n+k)$ of abelian varieties where the action of $K$ has
eigenspaces of dimension $n$ and $n+k$.

As before let $Y$ be a smooth ample divisor on $A$. The spaces $W^m$ are
Hodge substructures of level $m+k$ of $H^{g-m} (A ,\bQ)$ and also
$H^{g-m} (Y ,\bQ)$ via the embedding $H^{g-m } (A ,\bQ)\inj H^{g-m}
(Y ,\bQ)$. We have the
following refinement of Theorem \ref{thmGHT}

\begin{theorem}\label{weilhodge}
Suppose $m\leq n-1$.

If the Hodge conjecture holds for $W^m\subset H^{g-m} (A, \bQ)$, then
the Hodge conjecture holds for $W^{m+1}\subset H^{g-m-1}
(A, \bQ)$ and also for
$W^{m+1}\subset H^{g-m-1 }(Y ,\bQ)$.

If the Hodge conjecture is true for $W^{m+1}\subset H^{g-m-1 }(Y
,\bQ)$ and, for every subvariety $Z$ of pure codimension $n-m-1$ of $Y$,
the conjecture $GHC(C\times\tZ ,m+k+2, 1)$ holds, then the Hodge
conjecture holds for $W^m\subset H^{g-m} (A, \bQ)$.

In particular, using Lefschetz' Theorem in the case $k=m=0$, the Hodge
conjecture for $W^0\subset H^g (A,\bQ)$ is equivalent to that for
$W^1\subset H^{ g-1 } (A ,\bQ )$ and to that for $W^1\subset H^{ g-1 }
(Y ,\bQ )$.
\end{theorem} 

\begin{proof}
We first assume that the
Hodge conjecture holds for $W^m$ and prove it for $W^{m+1}$. There is then
an irreducible subvariety, say $Z_0$, of pure codimension $n-m$ of $A$ such
that $W^m$ is contained in the image of Gysin push-forward
\[
H^{ m+k } (\tZ_0 ,\bQ )\lra H^{g-m } (A ,\bQ).
\]
Let now $C$ be a smooth irreducible curve with a generically injective
map $C\ra A$ such that the image of $C$ generates $A$ as a group. Then
the pull-back map $H^1 (A ,\bQ )\ra H^1 (C ,\bQ )$ is injective. The
translates of $Z_0$ in $A$ by points of $C$ form the family of
subvarieties
\[
\begin{array}{ccc}
C\times Z_0 & \stackrel{q}{\lra} A \\
\downarrow^p & & \\
C & &
\end{array}
\]
whose image in $A$ we denote by $Z_1$. We choose $\tZ_1$ to be
$C\times\tZ_0$. Pontrjagin product and a
straightforward linear algebra computation now shows that the image of
\[
H^1 (A ,\bQ )\otimes H^{ m+k } (\tZ_0 ,\bQ )\subset H^1 ( C ,\bQ
)\otimes H^{ m+k } (\tZ_0 ,\bQ )\subset H^{ m+k+1 } (C\times\tZ_0 ,\bQ
)
\]
by Gysin push-forward is a Hodge substructure of $H^{ g-m-1 } (A ,\bQ
)$ containing $W^{m+1}$. So we obtain the Hodge conjecture for
$W^{m+1}\subset H^{ g-m-1 } (A,\bQ)$. The Hodge structure
$W^{m+1}\subset H^{ g-m-1 }(Y,\bQ )$ is now supported on $Z_1\cap Y$
which is a subvariety of pure codimension $n-m-1$ of $Y$ after replacing
$Z_1$ by a general translate if necessary to ensure that the
intersection of $Z_1$ with $Y$ is proper.

Assume now that the Hodge conjecture holds for $W^{m+1}\subset H^{ g-m-1
}(Y ,\bQ )$ and let $Z\subset Y$ be a subvariety of pure codimension
$m-n-1$ supporting $W^{m+1}$.
Consider the addition map $C\times Y\ra A$. The Hodge structure $H^1
(A ,\bQ )\otimes W^{m+1}$ is supported on the subvariety $C\times Z$ of
$C\times Y$. It is immediately seen that pull-back on cohomology from
$A$ induces an embedding
\[
W^m\inj H^1 (A,\bQ )\otimes W^{m+1}\inj H^1
(C ,\bQ)\otimes H^{m+k+1} (\tZ ,\bQ ).
\]
Now we use $GHC(C\times\tZ ,m+k+2, 1)$ to deduce that there is a
codimension $1$ subvariety $Z'$ of $C\times \tZ$ supporting $W^m$. The
pushforward $Z_0$ to $A$ of $Z'$ will then support $W^m$. This proves
the Hodge conjecture for $W^m\subset H^{ g-m } (A ,\bQ)$.

\end{proof}

It is well-known--and easy to prove--that the general Hodge conjecture
is equivalent to the
following statement
\begin{conjecture}
For every $\bQ$-Hodge substructure $V$ of $H^m (X ,\bQ)$ with level
$\leq m-2p$, there exists a nonsingular projective family of
subvarieties of pure dimension $g-m+p$ of $X$
\[
\begin{array}{ccc}
\cZ & \stackrel{q}{\lra} & X \\
\downarrow^r & & \\
S & &
\end{array}
\]
whose base is a nonsingular projective variety $S$ of dimension $m-2p$
such that the image of $H^{ m-2p } (S ,\bQ)$ by the Abel-Jacobi map
$q_* r^*$ of the family contains $V$.
\end{conjecture}

This formulation is a more geometric way of looking at the Hodge
conjecture and we use this formulation in the corollary below for
abelian varieties of Weil type $0$.

The corollary shows, for instance when $n=2$, that, if the Hodge
conjecture is true for $W^1\subset H^3 (Y ,\bQ)$, then it has a
solution given by a family of curves that can be embedded, fiber by
fiber, in a family of smooth surfaces in $Y$ such that the cohomology
class of each curve in its surface is linearly independent from the
restriction of the class of $Y$. So to solve the Hodge conjecture for
$W^1\subset H^3 (Y ,\bQ)$ or, equivalently, for $W^0\subset H^4 (A
,\bQ)$, we need to look for families of smooth surfaces in $Y$ with
Picard number $>1$. In higher dimensions we would look for families of
$n$-dimensional complete intersections in $Y$ that contain Weil divisors
homologically independent of the restriction of $Y$ and whose singular
loci are ordinary double loci of dimension $n-3$.

We say that $D$ has an ordinary double locus of dimension $n-3$ when,
either $n\geq 3$ and locally analytically $D$ is
isomorphic to the product of a smooth variety of dimension $n-3$ and a
threefold with an ordinary double point, or $n=2$ and $D$ is smooth.

\begin{theorem}
Suppose $A$ is a general abelian variety of Weil type and suppose that
the Hodge conjecture holds for $W^1\subset H^{2n-1} (Y ,\bQ)$. Then
there is a one-parameter family of subvarieties of dimension $n-1$ of
$Y$
\[
\begin{array}{ccc}
\cZ & \lra & Y \\
\downarrow & & \\
X & &
\end{array}
\]
whose Abel-Jacobi map is a solution to the Hodge conjecture for $W^1$
and satisfies the following property. For a general point $x\in X$,
any $N>> 0$, any choice of $n-1$ general divisors $D_1,\ldots , D_{
n-1}$ of the linear system $|\cI_{Z_x} (N Y) |_Y|$ the
singularities of the intersection $D := D_1\cap\ldots\cap D_{ n-1 }$
consist exactly of an ordinary double locus of dimension $n-3$ and the
homology class of $Z_x$ in $D$ is not a multiple of $[Y]_D$. Here
$Z_x$ is the fiber of $\cZ$ at $x$ and $\cI_{ Z_x }$ is the ideal
sheaf of $Z_x$ in $A$.
\end{theorem}

\begin{proof}

By Theorem \ref{weilhodge} the Hodge conjecture also holds for
$W^0$. Let $Z_0$ be a
subvariety of $A$ whose class is $\alpha +\lambda [Y]^n$ (for some
$\alpha\in W^0$, $\alpha\neq 0$ and $\lambda\in\bQ$) and let $C$ be a
smooth curve with a generically injective map $C\ra A$ whose image
generates $A$ as a group. Then, as in the proof of Theorem
\ref{weilhodge}, using Pontrjagin product and a straightforward linear
algebra argument, one can see that the family $\{Z_x := (Z_0
+x)\cap Y : x\in C\}$ gives an answer to the Hodge conjecture for a
summand of $[Y]^{ n-1 }\wedge H^1 (A ,\bQ )\oplus W^1\subset H^{ 2n-1
}( Y ,\bQ)$ which projects onto $W^1$.

By \cite{kleiman69} we can assume $Z_0$ to be smooth. For $N$
sufficiently large, the variety $Z_0$ is cut out scheme-theoretically
by divisors in the linear system $|\cI_{ Z_0 } (N Y)|$. In particular,
any $n-1$ general divisors $E_1 ,\ldots ,E_{ n-1 }$ in $|\cI_{ Z_0 }
(N Y )|$ intersect properly where $\cI_{ Z_0 }$ is the ideal sheaf of
$Z_0$ in $A$. Furthermore, for $N$ large enough,
\[
H^1 (\cI^2_{ Z_0 } (N Y )) =0,
\]
hence the restriction map
\[
H^0 (\cI_{ Z_0 } (N Y ))\lra H^0 \left(\frac{\cI_{ Z_0 }}{\cI^2_{ Z_0
}} (N Y )\right)
\]
is surjective so that any general $n$ sections of $\cI_{ Z_0 } (n Y )$
map to $n$ general sections of $\frac{\cI_{ Z_0 }}{\cI^2_{ Z_0 }} (n Y
)$. Let $s_1 ,\ldots , s_{ n-1 }$ be sections with divisors of zeros
equal to $E_1 ,\ldots , E_{ n-1 }$ respectively. As in the proof of
Theorem 4.2 in \cite{thomas05}, by our hypotheses of genericity the
singular locus of the intersection $E_1\cap\ldots\cap E_{ n-1 }$ is either
empty or is an ordinary double locus and is the locus where the images
of $s_1 ,\ldots , s_{ n-1 }$ in $H^0 (\frac{\cI_{ Z_0 }}{\cI^2_{ Z_0
}} (N Y ))$ fail to be independent. Again, for sufficiently large $N$
and because $\frac{\cI_{ Z_0 }}{\cI^2_{ Z_0 }} (N Y )$ has rank $n$
and the sections are general, this locus is either empty or has pure
dimension $n-2$. The intersection of $E_1\cap\ldots\cap E_{ n-1 }$
with a general translate of $Y$ will therefore either
have an ordinary double locus of dimension $n-3$ or be smooth. By
\cite{goreskymacpherson88} page 199, the pushforward map
\[
H_2 (E_1\cap\ldots\cap E_{ n-1 }\cap Y,\bQ )\lra H_2
(E_1\cap\ldots\cap E_{ n-1 },\bQ )
\]
is surjective. From this it follows that the homology class of $(Z_0
+x)\cap Y$ in $E_1\cap\ldots\cap E_{ n-1 }\cap Y$ is independent of
the homology class of $Y |_{ E_1\cap\ldots\cap E_{ n-1 }\cap Y
}$. The fact $H_2 (E_1\cap\ldots\cap E_{ n-1 }\cap Y ,\bQ)\neq H_2 (A
,\bQ)$ and the weak Lefschetz Theorem imply that $E_1\cap\ldots\cap
E_{ n-1 }\cap Y$ is {\em not} smooth if $n\geq 3$.

Now choosing our curve $C$ general and putting $D_i := E_i\cap Y$, the
family will have the desired property for
at least one choice of $D_1 ,\ldots , D_{ n-1}$ and therefore for all
sufficiently general choices of such divisors.

\end{proof}

\bibliographystyle{amsplain}

\providecommand{\bysame}{\leavevmode\hbox to3em{\hrulefill}\thinspace}
\providecommand{\MR}{\relax\ifhmode\unskip\space\fi MR }
\providecommand{\MRhref}[2]{%
  \href{http://www.ams.org/mathscinet-getitem?mr=#1}{#2}
}
\providecommand{\href}[2]{#2}

\end{document}